\documentclass[fleqn,12pt,twoside]{article}
\usepackage[headings]{espcrc1}

\readRCS
$Id: espcrc1.tex,v 1.2 2004/02/24 11:22:11 spepping Exp $
\newcommand{\AmS}{{\protect\the\textfont2
  A\kern-.1667em\lower.5ex\hbox{M}\kern-.125emS}}

\title{Vandermonde's quintic and  multiple decompositions of the number 1318}

\author{Jason A.C.~Gallas\address[IF]
    {Instituto de F\'\i sica,  
     Universidade Federal do Rio Grande do Sul,\\ 
     91501-970 Porto Alegre, Brazil 
}
 \thanks{Currently at: Rechnergest\"utzte~Physik~der~Werkstoffe, 
 ETH H\"onggerberg HIF~E12, Schafmattstrasse, 
 CH-8093~Zurich, Switzerland.}
}

\runtitle{Double decomposition of 1318}
\runauthor{J.A.C. Gallas}

\begin{document}

\maketitle

\begin{abstract}
This note records a curious numerical identity:
the number 1318, connected with Vandermonde's cyclotomic quintic,
may be decomposed in two distinct ways as a sum of 
products of pairs of numbers 
taken from the set \{$6, 16, 26, 41$\}, namely 
$1318 = 6\cdot41 + 16\cdot26 + 16\cdot41
      = 6\cdot16 + 6\cdot26  + 26\cdot41$.
Based on the existence of radical solutions of certain families of 
Abelian and generalized Abelian equations,
we conjecture the existence of an infinite number of analogous 
decompositions involving arbitrarily large sets of numbers.
\end{abstract}

\vspace{0.7truecm}


The object of this note is to call attention to a 
remarkable  multiple decomposition of the number
$1318  =2\cdot659$, namely
\[
\begin{array}{rrrrrrrrrrrrr} 
1318&=&6\cdot41 &+&16\cdot26 &+&16\cdot41  
    &=&246      &+&416       &+&656, \\
    &=&6\cdot16 &+&6\cdot26  &+&26\cdot41  
    &=&96       &+&156       &+&1066.
\end{array}
\]
These decompositions provide degenerate ways of representing $1318$ 
as sums of products of pairs of integers from the set
$S$ = \{$6,  16, 26, 41$\}.
Such integers satisfy 
\[ 6^2+16^2 +26^2+41^2 = 2649  \qquad\hbox{and}\qquad
    2649-1318=1331=11^3. 
\]
In addition, $2649-2 \cdot 1318 = 1331-1318 = 13$ is the next prime 
following $11$. 
All these numbers are needed repeatedly to simplify
long algebraic expressions leading to the quintic radicals defining
the roots of Vandermonde's celebrated cyclotomic quintic \cite{leb}
$V(x) = x^5 -x^4 -4x^3 + 3x^2 + 3x -1$, of discriminant $11^4$.

Multiplicity of decomposition arises from an elaborate cyclicity, 
or ambiguity, in combining algebraic expressions to reach intermediate
quintic radicals with norm $1318$:
\[
\begin{array}{rrrrrrrrrrrrr} 
  1318 &=& 6\cdot16 &+&26\,(6 + 41) &=& 96  &+& 26 \cdot 47 &=& 96  &+& 1222,\\
       &=& 6\cdot41 &+&16\,(26 + 41)&=& 246 &+& 16 \cdot 67 &=& 246 &+& 1072,\\
      &=& 16\cdot26 &+& 41\,(6 + 16)&=& 416 &+& 41\cdot 22  &=& 416 &+& 902,\\
      &=& 26\cdot41 &+& 6\,(16 + 26)&=& 1066&+& 6\cdot 42  &=& 1066 &+& 252.
\end{array}
\]
In fact, we find the existence of radical solutions of Vandermonde's 
quintic to be a direct consequence of the above ambiguity of factorizations
because they allow simplifications of huge algebraic  expressions 
to take place.
In particular, the decompositions above play a decisive role 
in obtaining  explicit analytical solutions for intricate nonlinear 
interdependency problems in discrete-time dynamical system \cite{jg1,jg2}.
This specific application along with the desire of obtaining multiple 
decompositions systematically motivates us to pose this problem to the
specialized audience.
We add that, based on explicit radical solutions obtained for certain 
families of Abelian and generalized Abelian equations \cite{jg2,jg3},
we conjecture the existence of an infinite number of analogous 
decompositions, involving arbitrarily large sets of numbers.

In algebra, the so-called Brahmagupta-Fibonacci identity 
implies that the product of two sums of two 
squares is itself a sum of two squares \cite{gj}:
\[
\begin{array}{rrr} 
  (a^2+b^2)(c^2+d^2) &=& (ac+bd)^2 + (ad-bc)^2,\cr
                     &=& (ac-bd)^2 + (ad+bc)^2.\cr 
\end{array}
\]
In other words, the set of all sums of two squares is closed under 
multiplication. 
Fermat's theorem on sums of two squares \cite{ot}
asserts that an odd prime number p can be expressed as 
$p = x^2 + y^2$ with integer $x$ and $y$ if and only if $p$ 
is congruent to $1\ (\hbox{mod } 4)$. 
Now, $1318  =2\cdot659$ and  $659 \hbox{ mod } 4 = 3$, thus
not expressible as a sum of two squares.
This seems to imply multiple decompositions not to be directly
related to the Brahmagupta-Fibonacci identity, a possible candidate
to crack the intricacies of multiple decompositions.

The algebraic cyclicity underlying the decompositions 
described here will be addressed elsewhere. 
They play an important role in the classification of ``shrimps'', 
i.e.~in the nucleation of stability islands in dissipative dynamical 
systems as discussed recently, 
for example, for discrete maps in the last paper of Lorenz \cite{lor},
or in continuous flows \cite{hub}.

\medskip

The author is supported by the Air Force Office of 
Scientific Research, Contract FA9550-07-1-0102, 
and by a  Research Fellowship from CNPq, Brazil.


\end{document}